\RequirePackage{fix-cm}
\documentclass[smallextended]{svjour3}       
\smartqed  
\usepackage{graphicx}
\usepackage{amsmath}
\newtheorem{lem}{Lemma}
\newtheorem{thm}{Theorem}

\newtheorem{rem}{Remark}
\newtheorem{defi}{Definition}

\newtheorem{excer}{Exercise}

\begin{document}

\title{Separating  circles on the sphere by polygonal tilings}


\author{Andr\'as Bezdek }


\institute{Work of author was supported by NKFIH grant KKP-133864\\
 A. Bezdek \at
              MTA R\'enyi Institute, Budapest, Hungary\\
              Department of Mathematics and Statistics, Auburn University, Auburn, AL, USA\\
              \email{bezdean@auburn.edu}           
                     }

\date{Received: date / Accepted: date}

\maketitle

\begin{abstract} We say that a tiling separates discs of a packing in the Euclidean plane, if each tile contains exactly one member of the packing. It is a known elementary geometric problem to show  that for each locally finite packing of circular discs, there exists  a separating tiling with convex polygons. In this paper we show that this separating property remains true for circle packings on the sphere and in the hyperbolic plane. Moreover, we show that in the Euclidean plane circles are the only convex discs, whose packings with similar copies can be  always separated by polygonal tilings. The analogous statement  is not true on the sphere and it is not known in the hyperbolic plane.

\keywords{circle packing \and polygonal tiling \and Voronoi diagram}
\subclass{52C15 \and 52A55}
\end{abstract}

\section{Introduction and outline of the paper}

Cell partitions of the plane play an important role in solving optimization problems in the area of discrete geometry (see Fejes Toth \cite{FejesToth} and Pach \cite {Pach}. Suppose we have a locally finite  family of  points or circular discs or general  convex discs in the plane.  If no two of the discs have an interior point in common, then the family is said to form a {\it packing}. The phrase {\it circle packing} will refer to a packing of circluar discs.  The  {\it Voronoi diagram} associated with a packing is a partition of a plane into cells, where each cell consists of all points of the plane closer to a particular member of the packing,  than to any other. Voronoi diagrams are also called {\it nearest point partitions}. In the simplest case, when a point set is given, the cells are convex polygons. Note that  some of the cells could be unbounded (half planes, strips, angular sectors), thus for the purpose of a this  paper, we say that {\it convex polygons} are convex regions bounded by a finite number of segments, lines or half lines.

\begin{defi} For a given packing of convex discs, a cell partition of the plane is  called a {\it separating polygonal tiling}, if
 each cell is a convex polygon in the above extended sense and contains exactly one packing element.
\end{defi}

\begin{figure}[h]
\begin{center}
\includegraphics[scale=.72]{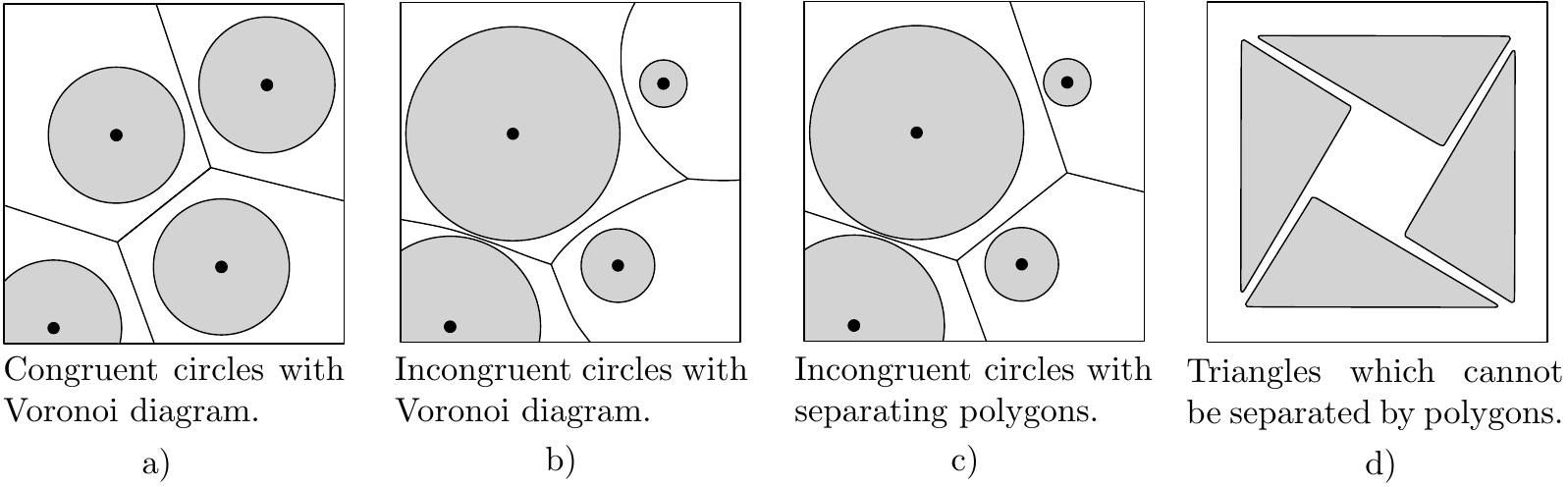}
\caption{Separating circles in the Euclidean plane}
\label{vor}
\end{center}
\end{figure}

If the circles in the packing are congruent, then the Voronoi diagram of the circles is the same as the Voronoi diagram of the centers of the circles, therefore it is a separating polygonal tiling (Figure \ref{vor}a). If the circles are not necessarily congruent, then the Voronoi diagrams  might have cells with curved boundaries (Figure \ref{vor}b), thus it is not necessarily a separating polygonal tiling. Moreover, the Voronoi diagram determined by the centers of the circles might have cells, which only partially contain circles.

In this paper, we study the existence of separating polygonal tilings for various packings. Our work was inspired by the following problem which appeared in a Russian math olympiad problem collection \cite{skljarszkij},

\begin{excer}\label{russian} Suppose we have finite many disjoint circular discs in a square. Prove that the square can be divided into convex polygons, so that each polygon contains exactly one disc (Figure \ref{vor}c ).
\end{excer}

Section 2 contains an outline of the solution of Exersice \ref{russian}. The proof presented in Section 2 shows also  that  spheres of  any locally finite sphere packing in $R^n$ can be separated by a polyhedral tiling. Figure \ref{vor}d is also included in \cite{skljarszkij}  to let readers realise that they are packings of convex discs, whose members cannot be separated by polygonal tilings. This paper starts with two additional generalizations of Exercise \ref{russian}.

\bigskip

\begin{figure}[h]
\begin{center}
\includegraphics[scale=.57]{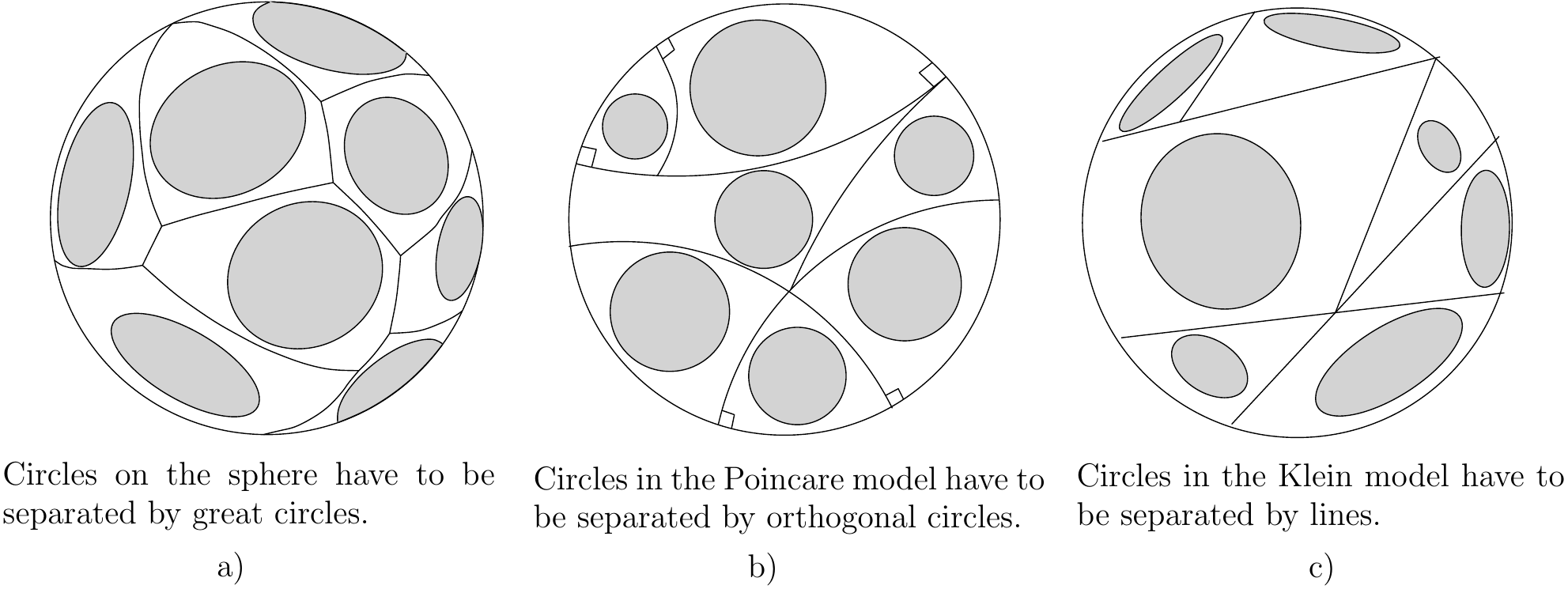}
\caption{Polygonal separation in the sphere and in the hyperbolic plane}
\label{voron2}
\end{center}
\end{figure}

The first real generalizations concerns the separability of circle packings on the sphere. Arcs of the great circles are the geodesics on the sphere, thus we want to  partition the sphere into convex regions  bounded by arcs of great circles (Figure \ref{voron2}a). Such tilings on the sphere will be called polygonal tilings. This automatically means that we consider only packings of circular discs which are  smaller than a hemisphere. A very natural idea would be to use a stereographic projection, say from the north pole to a plane tangent to the sphere at the south pole. This projection maps circles of the sphere to planar circles. In view of Exercise \ref{russian}, the planar circles can be separated by a polygonal tiling.   Pulling back the boundary segments of this tiling to the sphere, one would get a separating tiling on the sphere. The problem is that although the boundary of such tiles consists of arcs of circles, they are not necessarily arcs of great circles. In Section 3 a different approach is presented to prove

\begin{thm}\label{spheresep} For every  finite family of  disjoint circular discs of radii $\leq \frac \pi 2$ on a sphere there is a separating spherical polygonal tiling.
\end{thm}

The second generalization concerns circle packings in the hyperbolic plane. In the hyperbolic plane all the needed concepts (lines, circles, convexity, polygons) exist, thus one can ask the polygonal separation problem exactly in the same way as in the Euclidean plane. A natural approach would be to get a better understanding of the properties of the polygonal tiles either in the Poincare disc model or in the Klein disc model and use  Euclidean arguments to show the existence of separating polygonal tilings. An advantage of the Poincare disc model is that it is conformal (circles and angles are not distorted), a disadvantage is that lines of the geometry are circular arcs orthogonal to the boundary circle of the disc (Figure \ref{voron2}b). In the Klein disc model lines of the hyperbolic plane are chords of the disc.  Circles in the Klein disc model become ellipses increasingly flattened as they are nearer to the edge of the model (Figure \ref{voron2}c).  Thus in the Klein disc model while  lines are simple to study for separation, circles are not. In Section 4 a different approach is presented to prove

\begin{thm} \label{hyperbolicsep} For every  locally finite family of  disjoint circular discs in the hyperbolic plane there is a separating polygonal tiling.
\end{thm}

Finally, we turn to the  question of characterizing circles with the separating property. This part of the paper focuses to convex discs the Euclidean plane,

\begin{defi}We say that a convex disc $C$ has the separating property, if every packing of similar copies of $C$ can be separated by a polygonal tiling.
\end{defi}

We define 'caps' of convex discs, in particular 'isosceles caps' as they play an important role in such discussion.

\begin{defi} Let $C$ be a convex disc and let $S$ be a point outside of $C$. Assume that the two tangent lines of $C$ through $S$ touch $C$ at points $P$ and $Q$. Cap $\widehat{PSQ}$ is the union of segments $PS$ and $SQ$. Cap $\widehat{PSQ}$ is called {\it isosceles}, if the segments $PS$ and $SQ$ have equal lengths.
\end{defi}

Section 5 contains two theorems concerning isosceles caps of a convex discs in the Euclidean plane.  These theorems have analogous versions on the sphere and in the hyperbolic plane, but we do not include them in Section 5, as they are not going to be applied in Section 6.

\begin{thm} \label{circle lemma} Circles are the only convex discs in the Euclidean plane, whose caps are all isosceles.
\end{thm}

\begin{thm} \label{isoscelescap} Let $0 < \alpha < \pi$. Every convex disc $C$ in the Euclidean plane  has an isosceles cap of angle $\alpha$.
\end{thm}

The above two theorems will be applied in Section 6, where  the following theorem will be proved,

\begin{thm}\label{characterization} Let $D$ be a noncircular disc in the Euclidean plane. Finite many disjoint similar copies of disc $D$ can be arranged in the plane so that they cannot be separated by a polygonal tiling.
\end{thm}

Note that  similarity does not exists neither on the sphere nor in the hyperbolic plane, so one would need to restrict the question of separation for packings of congruent copies of discs in those two spaces. Even under this restriction, Theorem \ref{characterization} is not true on the sphere. Just take any convex disc inside of a hemisphere, with an area greater than $\frac 23$ of the surface area of the sphere. Obviously, there is room only for at most two copies of $C$ on the sphere, and two copies can always be separated by a great circle. In the hyperbolic plane it is not known if circles are characterised with the separation property.

\section{Separation in the Euclidean plane, proof of Exercise \ref{russian}}

In this section, we include the outline of a known solution of Exercise \ref{russian}, so that the reader could see how the generalizations are related to this original solution.

\begin{defi} The power of a point $A$ with respect to a circle $C$ centered at $O$ of radius $r$ is defined by
 $p(A,C)=AO^{2}-r^{2}$.
 \end{defi}

By this definition, points inside the circle have negative power, points outside of the circle  have positive power, and points on the circle have zero power. For external points, the power equals the square of the length of a tangent from the point to the circle.

\begin{lem} \label{isopower} Let $C_1$ and $C_2$ be two disjoint circles. The points which have equal power with respect to $C_1$ and $C_2$ form a line. This line is perpendicular to the line containing the centers of the circles and separates the circles.
\end{lem}

\begin{figure}
\begin{center}
\includegraphics[scale=.5]{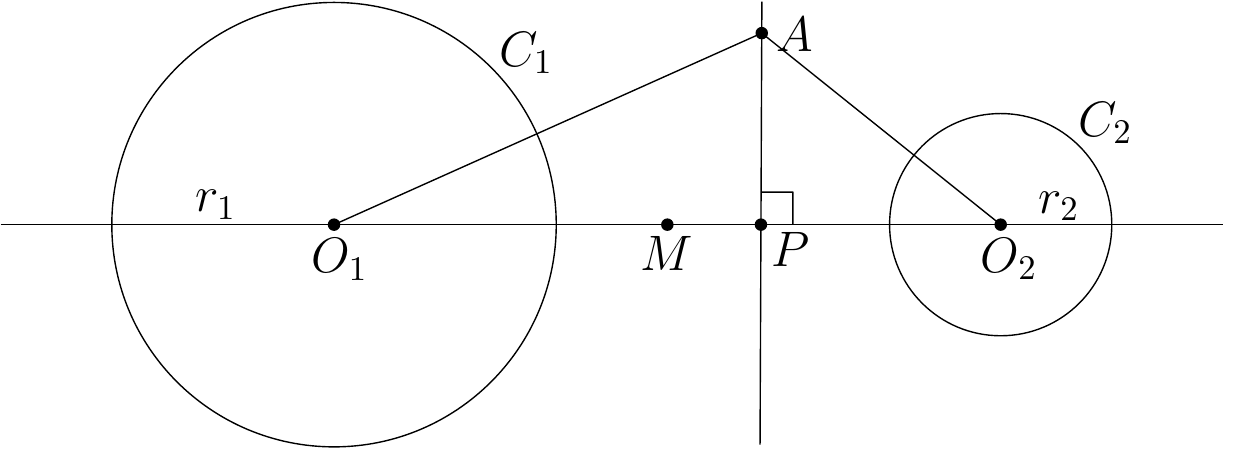}
\caption{Points having equal power with respect to two circles}
\label{power}
\end{center}
\end{figure}

\noindent {\bf Proof of Lemma \ref{isopower}.} Let $O_1$ and $O_2$ be the centers the circles $C_1$ and $C_2$  and let  $r_1 \geq r_2$ be the radii of the circles $C_1$ and $C_2$ (Figure \ref{power}). Applying Pythagorean theorem for the right triangles $O_1PA$ and $O_2PA$, we have that if a point $A$ has equal powers with respect to the two given circles, then so does the perpendicular projection $P$ of the point $A$ on the line connecting $O_2$ and $O_1$. Since $r_1 \geq r_2$, $P$ belongs to the ray $MO_2$ where $M$ is the midpoint of the segment $O_2O_1$. It is easy to verify that the power difference of point $P$ with respect to the two circles is
$p(P, C_1) -p(P,C_2) = (PO_1-PO_2)(PO_1+PO_2) - r_1^2 +r_2^2 = 2  MP \cdot  O_1O_2 -(r_1^2 -r_2^2)$. This equation implies that there is exactly one point $P$ on the ray $MO_2$ for which the power difference is zero. Thus, the points which have equal power with respect to $C_1$ and $C_2$ form a line, which  is perpendicular to the line containing the centers of the circles. In fact, on that side of this line which contains $O_1$, lie those points whose power is grater with respect to circle $C_1$, while on the other side the situation is exactly the opposite, which implies the needed separation of circles $C_1$ and $C_2$. \qed

Similarly to the Voronoi partition,  Lemma \ref{isopower} implies that the least power partition is a separating polygonal tiling.

\section{Separation on the sphere, proof of Theorem \ref{spheresep}}

In this section we prove Theorem \ref{spheresep}, which is formally stated in the introduction. In plain terms, we are going to construct a separating polygonal tiling for every finite circle packing on the sphere.

For the purpose of the construction we introduce the {\it spherical potential} of a point with respect to a given circle. First we review some standard notations: distance $\widehat{PQ}$ between two points $P,Q$ on the sphere is measured by the central angle of the shorter great circular arc connecting the two points. Similarly, radius of a circle is measured by the central angle corresponding to the arc representing the radius of the circle.

\begin{defi} The spherical potential of a point $A$ with respect to a circle $C$ centered at $O$ of radius $r$ is defined by $p(A,C) = \frac {\cos \widehat{AO}} {\cos r}$.
\end{defi}

The spherical potential was introduced specifically for this proof. The potential $p(A,C)$ depends on the distance $\widehat{AO}$ through cos(.).  Further the point is on the sphere from the center, less potential it has.

We prove the following  two lemmas.

\begin{lem}\label{equipotential1}   Assume that two disjoint circular discs are given on a sphere $S$ and their centers lie on a great circle $C$. Let $A$ be a random point on $S$ different from the centers of the hemispheres bounded by $C$. Let $P$ be the perpendicular projection of $A$ onto the great circle $C$. Then $A$ is an equipotential point with respect to the two given  circles, if and only if the projection $P$ is also an equipotential point.
\end{lem}

\begin{figure}[h]
\begin{center}
\includegraphics[scale=.6]{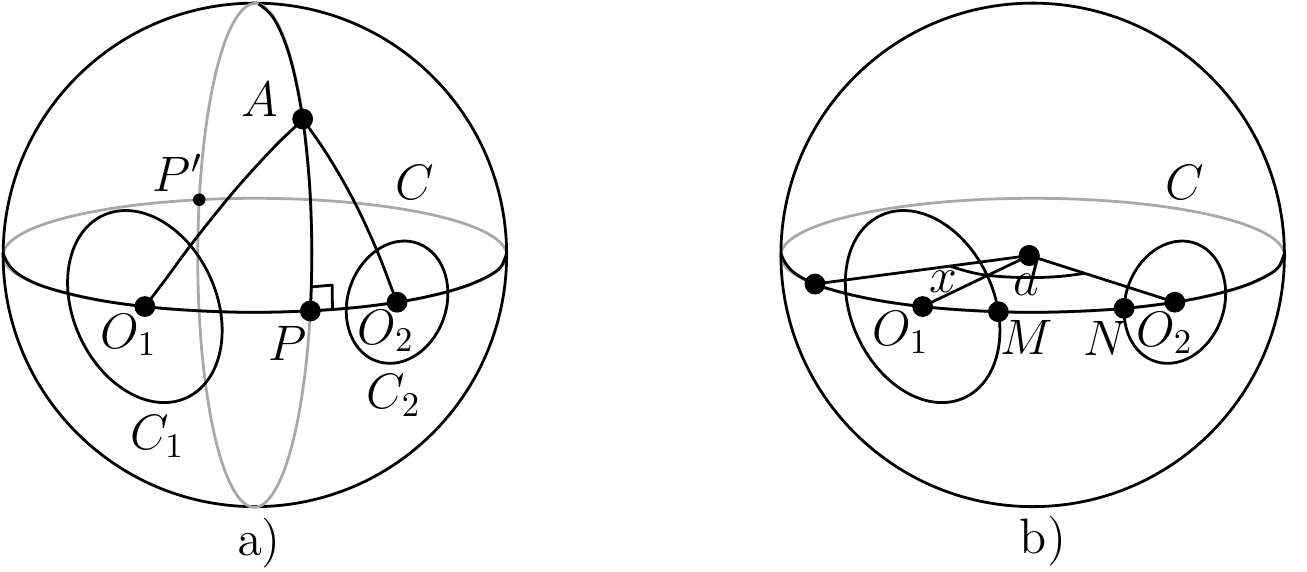}
\caption{Equipotential points are on great circle}
\label{equiline}
\end{center}
\end{figure}

\noindent {\bf Proof of Lemma \ref{equipotential1}} Let $C_1, C_2$ be two given circles with centers $O_1, O_2$ and radii $r_1, r_2$ (Figure \ref{equiline}a).   Note that the  centers of the two hemispheres bounded by $C$ have $0$ potentials with respect to both circles. If $A$ is different from these two 'poles' then project $A$ onto the great circle $C$ passing through  $O_1$ and $O_2$. Let $P$ be the projection of point $A$ on $C$ are the projections.
Since $A$ is an equipotential point, we have

$$ \frac {\cos \widehat{AO_1}} {\cos r_1} = \frac {\cos \widehat{AO_2}} {\cos r_2}$$

According to the spherical law of cosine, in any right triangle with hypotenuse $c$ and legs $a,b$ we have $\cos a \cos b = cos c$. Thus, the numerators in the previous equation can be replaced with appropriate products of cos's:

$$ \frac {\cos \widehat{AP} \cdot \cos \widehat{PO_1}} {\cos r_1} = \frac {\cos \widehat{AP} \cdot \cos \widehat{PO_2}} {\cos r_2}$$

After cancelling with $\cos \widehat{AP}$ we get that P is also an equipotential point.
\qed

\begin{lem} \label{equipotential2} Assume two disjoint  circular discs are given on a sphere. Then on the great circle  passing through the centers of them there ar exactly two points which are equipotential with respect to the two circles. Moreover, these points are antipodal and they separate the circles.
\end{lem}

\noindent {\bf Proof of Lemma \ref{equipotential2}}. {\bf Existence part od the proof:} Let $C_1, C_2$ be two given circles with centers $O_1, O_2$ and radii $r_1, r_2$ and let $C$ be the great circle passing through $O_1$ and $O_2$ (Figure \ref{equiline}b). Let the shorter great circular arc of $C$ intersect the circles $C_1$ and $C_2$ at points $M$ and $N$. By definition of the potential, we have that

$$p(M, C_1) = 1 \quad  \text{and} \quad p(M, C_2) < 1,$$

$$p(N, C_1) < 1 \quad \text{and} \quad p(N, C_2) = 1$$.

Moving a point P from M to N along  the arc MN, the inequality changes from $p(X, C_1) > P(X, C_2)$ to $p(X, C_1) < P(X,C_2)$. A simple continuity argument implies that  there is at least one  equipotential point on the arc MN.

\bigskip

\noindent {\bf Uniqueness part of the proof:} Let $d$ be the central angle of the shorter great circular arc $\widehat{PO_2}$ (Figure \ref{equiline}b). Let us move a point $P$ along the great circle $C$ starting at $O_1$ along the opposite direction of $O_2$. Let us identify the position of $P$ with the central angle  $x = \widehat{PO_1}$. To find all possible equipotential points on $C$ we need to solve the following equation for $x$ over the interval $[ 0, \pi ]$:

$$ \frac {cos x}{cos r_1} = \frac {cos (x+d)}{cos r_2},$$

which is the same as

$$ \frac {cos r_2}{cos r_1} = \frac {cos x \cdot \cos d - \sin x \cdot \sin d}{cos x},$$

Solving this for x, we have that

$$\tan x = \frac {\cos r_1}{\sin d \cdot \cos r_2} = constant.$$

Since this constant is not equal to $0$, we have a unique solution for $x$ on the half circle, i.e. on the interval $[ 0, \pi )$. Note that the antipodal point of this unique $P$ is also an equipotential point  with respect to the circles $C_1$ and $C_2$. \qed

Uniqueness of the perpendicular projections mean that all equipotential points lie on a great circle perpendicular to the great circle $C$.  Moreover, by monotonicity of the potential it turns out that this great circle is separating $C_1$ and $C_2$. Similarly to the Voronoi partition, we have now  that the least potential partition is a separating polygonal tiling, which completes the proof of Theorem \ref{spheresep}.

\section{Separation in the hyperbolic plane, the proof of Theorem \ref{hyperbolicsep}}

It turns out that after appropriately redefining the potential of a point $A$ with respect to a  circle $C$ of center $O$ and radius $r$, the entire proof of the spherical separation problem goes through word by word. In fact, the potential $p(A,C) = \frac {\cosh \widehat{AO}} {\cosh r}$ does the job.

\section{Caps of convex discs, proof of  Theorem \ref{circle lemma}}

\noindent Theorem \ref{circle lemma} claims that  circles in the Euclidean plan are the only convex discs whose caps are all isosceles.

\noindent {\bf 1st proof  of Theorem \ref{circle lemma}:}` Let $C$ be a convex disc whose caps are all isosceles. We may assume that the $C$ is strictly convex, otherwise $C$ has a tangent line with several contact points, permitting non isosceles caps. Let us start with recalling a simple elementary exercise, where one has to find a point on each side of a given triangle so that each vertex is at equal distances from the selected points of its adjacent sides. Note that the triplet of tangency points of the incircle of the given triangle is a solution. Moreover, this solution is unique. Indeed, if the lengths of the side partitions clockwise are $x, x, y, y, z, z$, then decreasing one of the lengths $\{x,y,z\}$ forces the other two lengths  to increase, while their sum should stay constant.

\begin{figure}[h]
\begin{center}
\includegraphics[scale=.6]{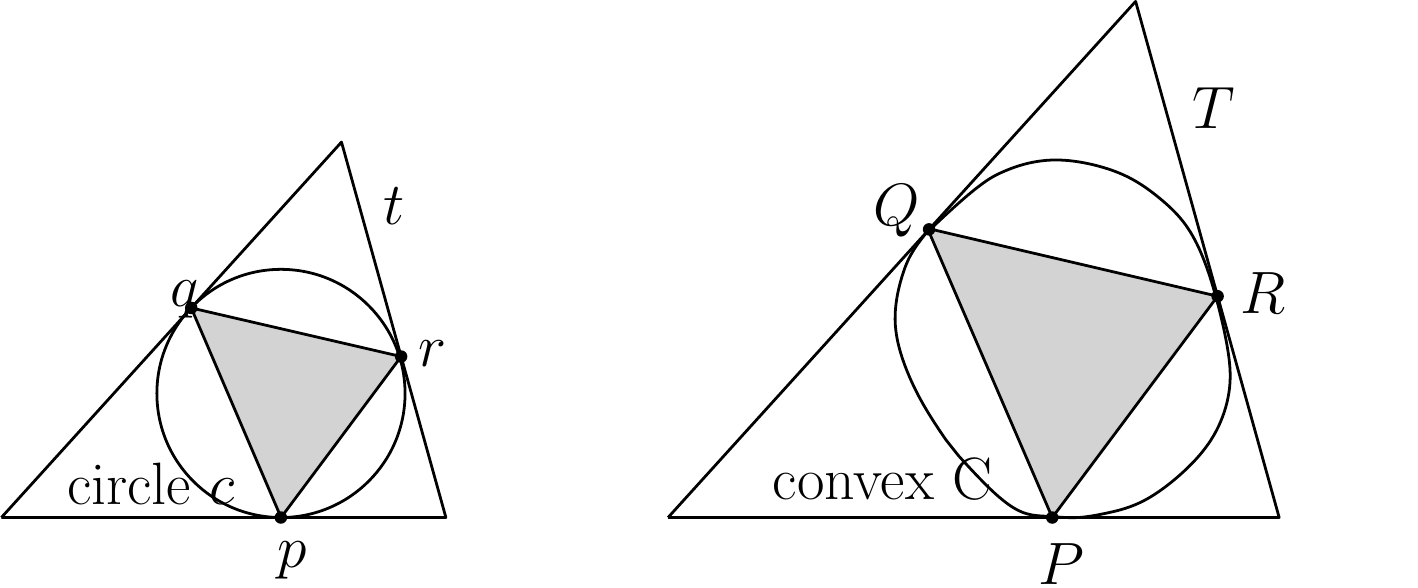}
\caption{Comparing a circle $c$ to a convex $C$ with isosceles caps}
\label{definition}
\end{center}
\end{figure}

Returning to $C$, for comparison put a circle $c$ next to $C$ (Figure  \ref{definition}). Let $p$ and $P$ be the lowest points of $c$ and $C$. Consider a random circumscribed triangle $t$ with contact points $p, q, r$ of circle $c$. Let $T$ be the circumscribed triangle of $C$, whose sides are parallel to the sides of triangle $t$. Denote the corresponding tangency points of $T$ and $C$ by $P, Q, R$. In view of the above elementary exercise, the triangles $pqr$ and $PQR$ are not only similar, but they are in parallel positions. Let us now fix two of the contact points $\{p, q, r\}$ (note that there are three ways to do this) and change triangle $t$ by  running the third contact point  along the arc of circle $c$. The later observation implies that $R$ runs along a circular arc. Moreover, the three arcs corresponding to the three choices of the two contact points of $\{p,q,r\}$, form a complete circle. \qed

Note that the above proof uses similarity, thus it does not work on the sphere and in the hyperbolic plane. Although we are not going to use the spherical and hyperbolic versions of Theorem \ref{circle lemma}, we present here a second Euclidean proof, which can be modified for those two spaces.

\bigskip

\begin{figure}[h]
\begin{center}
\includegraphics[scale=.4]{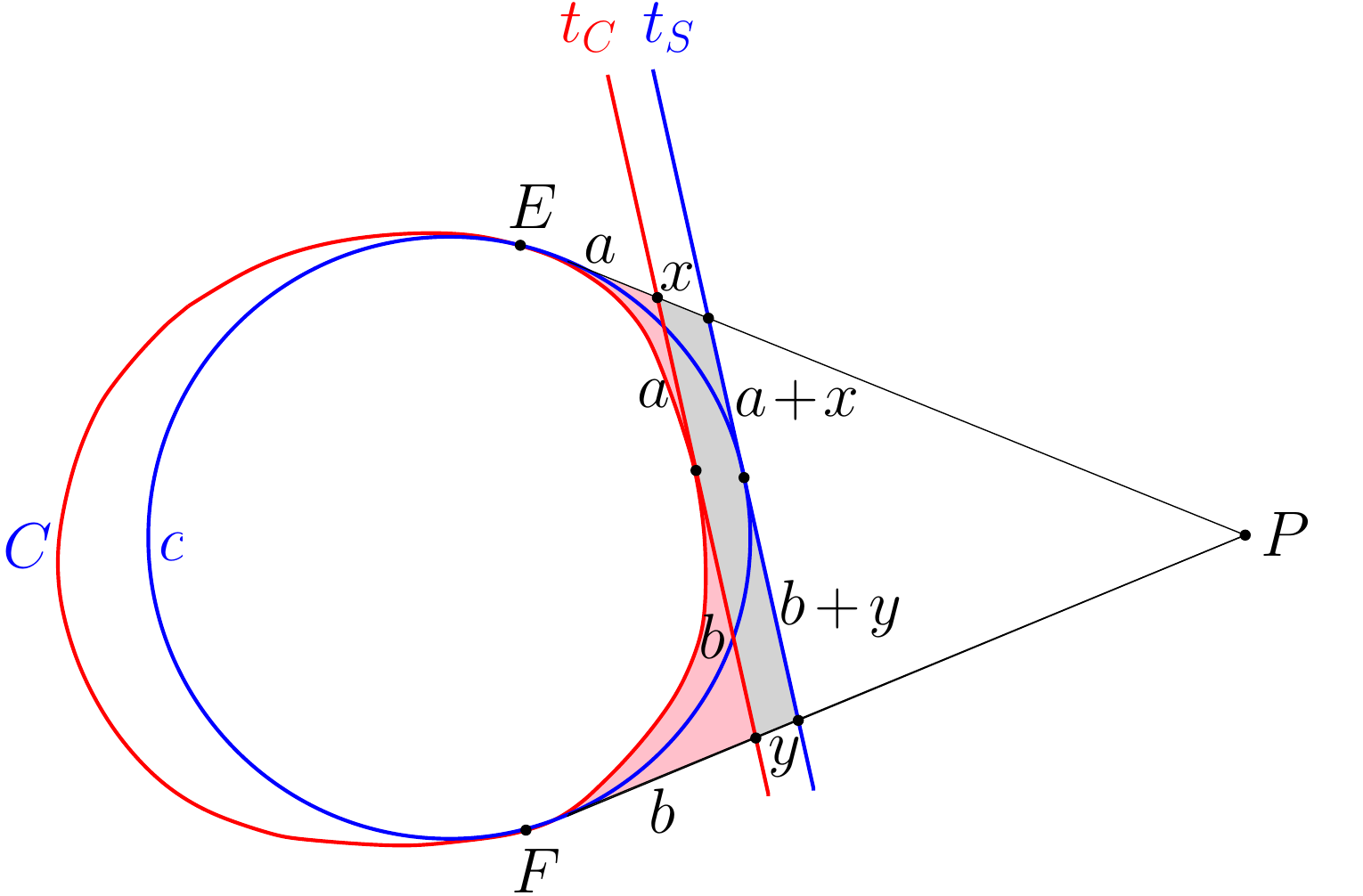}
\caption{Comparing a circle $c$ to a convex $C$ with isosceles caps}
\label{circle2}
\end{center}
\end{figure}

\noindent {\bf 2nd proof of Theorem \ref{circle lemma}:} We may assume that the convex disc $C$ is strictly convex. Otherwise, $C$ has a tangent line with several contact points, permitting non isosceles caps. Let again $C$ be a convex disc whose caps are all isosceles. Consider a cap of $C$, with vertex $P$ and sides $EP$ and $FP$.  Since all caps are isosceles, there is a unique circle $c$, which  is tangent to  the lines $EP$ and $FP$. Indirect assume  $C$ is different from circle $c$.   $C$ must have  a tangent line $t_C$, which  is not a tangent line of circle $c$ (Figure  \ref{circle2}). Shift $t_C$, without passing through $E$ and $F$, until it becomes  tangent to circle $c$. Label this tangent line by $t_c$. Note, that with this carefully phrased instruction we cover all subcases  i) $t_C$ intersects circle $c$ or not, ii) $t_C$ separates $E$ from $P$ or not.  Let $R$ be  the quadrilateral  bounded by the lines $t_c, t_C, EP$ and $FP$. Let $x$ and $y$ denote the lengths of the opposite sides of $R$ along the segments $EP$ and $FP$. Let  $a$ and $b$ be the lengths of the sides  of those two  caps,  whose vertices are on the lines $EP$ and $FP$ and are closest to $E$ and $F$. Using these notations the sides of quadrilateral $R$ have lengths $x, (a+b), y,  (a+x + b+y)$. Since the sum of the first three lengths is equal to the  fourth length, quadrilateral $R$ is degenerated, thus $t_c$ coincides with $t_C$ , a contradiction.  \qed

\begin{figure}[h]
\begin{center}
\includegraphics[scale=.4]{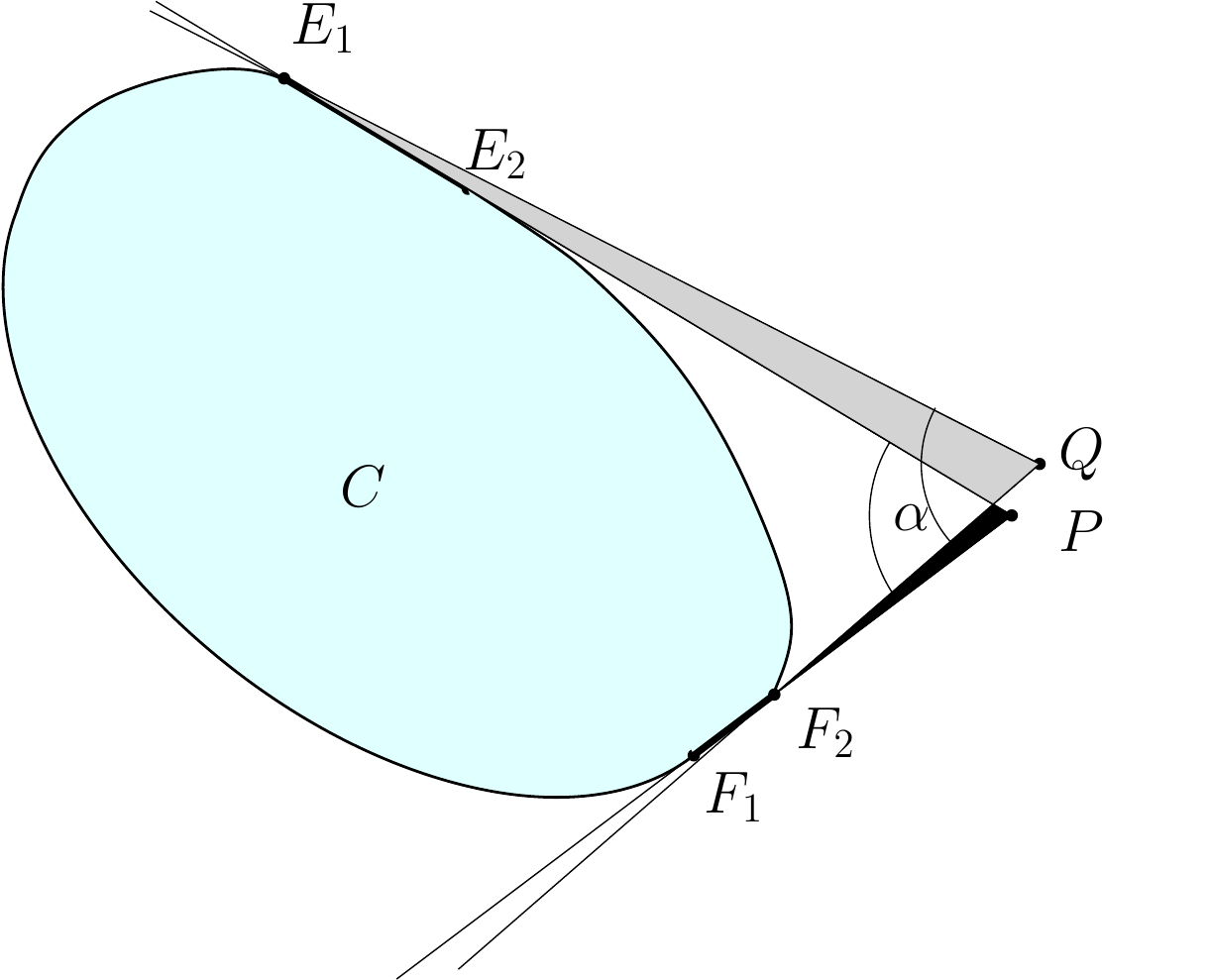}
\caption{There exists isosceles cap of given angle}
\label{iso-cap}
\end{center}
\end{figure}

\noindent {\bf Proof of Theorem \ref{isoscelescap}}.  Theorem \ref{isoscelescap} claims that for $0 < \alpha < \pi$, every convex disc $C$ in the Euclidean plane has an isosceles cap of angle $\alpha$. Indirect assume that there is an $\alpha$ so that all caps of angle $\alpha$ of disc $C$ are scalene. A simple compactness and continuity argument shows that among caps of angle $\alpha$ there is one with maximum area. Since we do not assume strict convexity we have to allow the possibility that the sidelines of this cap share a segment $E_1 E_2,  F_1 F_2$ with the boundary of $C$. Without loss of generality, we may assume that $E_1P \geq E_2P >  F_1P  \geq  F_2P$. (Figure  \ref{iso-cap}). In short side $EP$ is longer then side $FP$. Let us roll the support line $E_1E_2P$ on $C$ away from cap $EPF$ and roll the support line $E_1E_2P$ on $C$ into the cap $EPF$.  If the support lines are  rotated by the same small angle then they determine a new cap of angle $\alpha$. On  Figure  \ref{iso-cap} the area of the original and the rotated caps can be compared.  The gray  shaded area represents area of gain, the black shaded are represents the are of loss. Since $E_2P >  F_1P$ one can easily see that the area of gain is larger, than the area of loss, which contradicts the maximum area choice. \qed

\section{Non circular discs admit non separable packings}

\begin{figure}[h]
\begin{center}
\includegraphics[scale=.6]{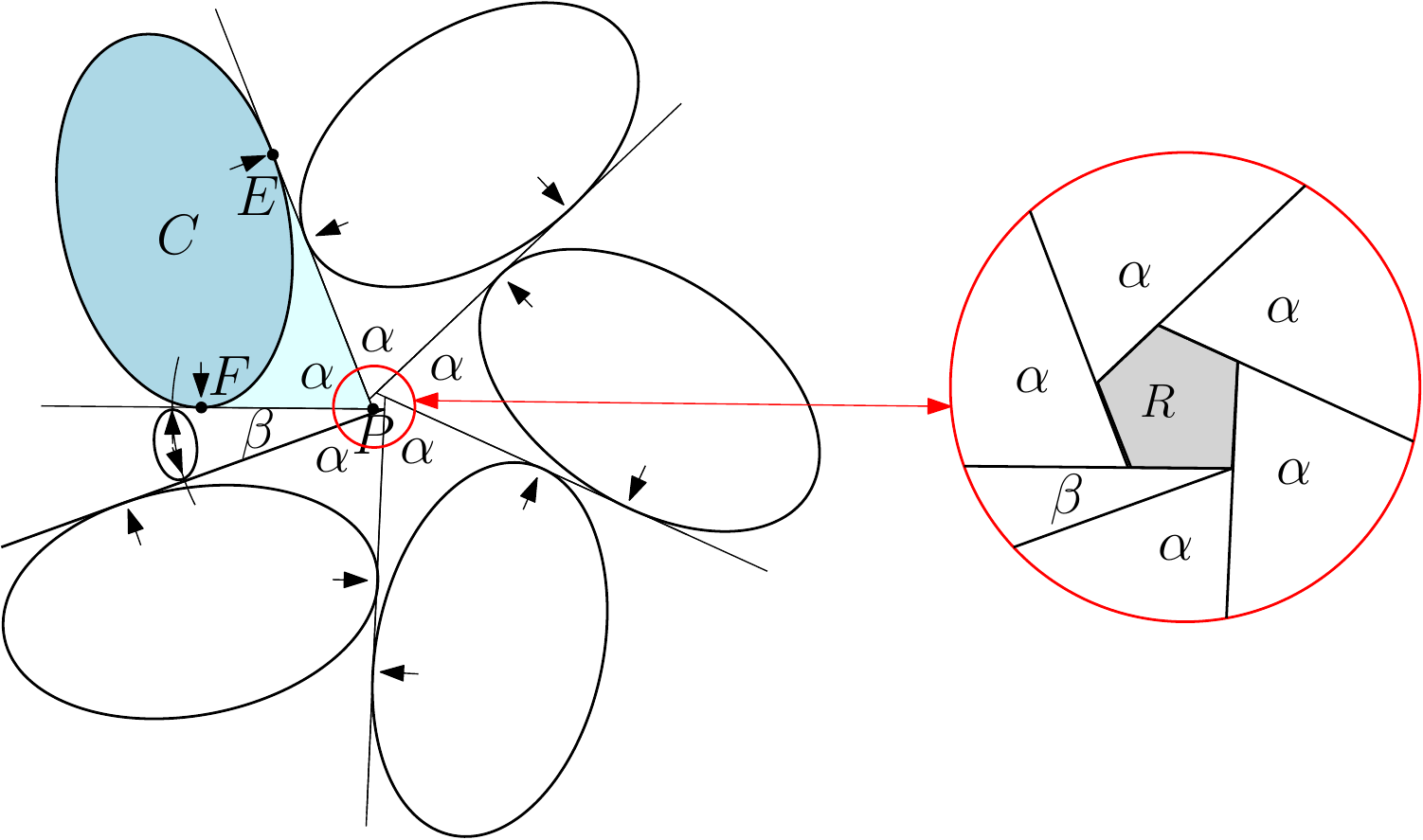}
\caption{Nonseparable packing of convex discs}
\label{counter example}
\end{center}
\end{figure}

In this section we prove Theorem \ref{characterization}, which is formally stated in the introduction. With other words, we will give instructions how to arrange finite many similar copies of a given non circular discs, which already form a nonseparable packing.  Assume the convex disc $C$ is not a circle. In view of Theorem \ref{circle lemma}, $C$ has a non isosceles cap $EPF$ with $EP > FP$. For future reference chose $\epsilon >0$ to be smaller than $EP-FP$. Let $\alpha$ be the angle of the cap EFP. Let $n$ be the largest  integer so that $n\alpha \leq 2 \pi$ and let $\beta = 2\pi - n\alpha$. For our construction we need to start with an $(n+1)$ sided polygon $R$, whose outside angles are all equal to $\alpha$ except one which is equal to $\beta$ and whose sides are all shorter than $\epsilon$.  Figure \ref{counter example} shows an arrangement of $n$ congruent copies of $C$ so that the attached $\alpha$ caps  together with an empty $\beta$ angular sector tile the neighbourhood of polygon $R$. According to Theorem \ref{isoscelescap} disc $C$ has an isosceles $\beta$ cap.  Finally, we place a scaled copy of $C$ with the isosceles $\beta$ cap fitting in the $\beta$ angular sector. By scaling we can ensure that the side lengths of the $\beta$ cap are between  $EP-\epsilon$ and $FP$. The choice of $\epsilon$ and the clockwise rotational condition guaranties that no separating line between consecutive copies of $C$ can not intersect the polygon $R$, thus separating tiling cannot exists. \qed

\end{document}